\documentclass[12pt,draftcls,onecolumn]{IEEEtran}

\usepackage{amsfonts}
\usepackage{amsmath}
\usepackage{amsfonts}
\usepackage{amssymb}
\usepackage{psfrag}
\usepackage{graphicx}
\usepackage{multirow}

\def\diag{\mbox{diag}}

\def\ra{\rightarrow}

\def\ra{\rightarrow}

\def\np{\newpage}

\def\matlab{MATLAB$^{\textrm{\tiny{\textregistered}}}$}
\def\trademark{\textrm{\tiny{\texttrademark}}}
\def\registered{\textrm{\tiny{\textregistered}}}

\def\bd{\begin{displaymath}}
\def\ed{\end{displaymath}}
\def\bi{\begin{itemize}}
\def\ei{\end{itemize}}
\def\bn{\begin{enumerate}}
\def\en{\end{enumerate}}
\def\bq{\begin{eqnarray}}
\def\eq{\end{eqnarray}}
\def\bqn{\begin{eqnarray}}
\def\eqn{\end{eqnarray}}
\def\bqta{\begin{eqtarraya}}
\def\eqta{\end{eqtarraya}}
\def\bqtb{\begin{eqtarrayb}}
\def\eqtb{\end{eqtarrayb}}
\def\bqtc{\begin{eqtarrayc}}
\def\eqtc{\end{eqtarrayc}}
\def\be{\begin{equation}}
\def\ee{\end{equation}}
\def\bea{\begin{eqnarray}}
\def\eea{\end{eqnarray}}
\def\beann{\begin{eqnarray*}}
\def\eeann{\end{eqnarray*}}

\newcommand{\real}{{\mathbb{R}}}

\newcommand{\rank}{{\rm rank}}

\newcommand{\ind}{{\rm index}}

\newcommand{\byers}{{\textit{byersnash}}}

\newcommand{\Span}{{\textit{span}}}
\newcommand{\sylv}{{\textit{sylvplace}}}
\newcommand{\place}{{\textit{place}}}
\newcommand{\rfbt}{{\textit{rfbt}}}
\newcommand{\rob}{{\textit{robpole}}}

\def\undp{{\underline \pi}}
\def\overp{{\overline \pi}}
\def\bmat{\left[ \begin{array}}
\def\emat{\end{array} \right]}

\def\bsmat{\left[ \begin{smallmatrix}}
\def\esmat{\end{smallmatrix} \right]}

\def\im{{\rm im}\ }

\newtheorem{prop}{Proposition}[section]

\begin{document}

\title {Robust pole placement with Moore's algorithm}
\author{Robert Schmid,  Amit Pandey  and Thang Nguyen  \thanks{Robert Schmid is  with  the  Department of Electrical and Electronic Engineering, University of Melbourne.   Amit Pandey is with  the  Department of Electrical and Computer Engineering, University of California at San Diego. Thang Nguyen is with  the Department of Engineering, University of Leicester.
email:
rschmid@unimelb.edu.au,  amitpandey.ent@gmail.com, tn57@le.ac.uk. An  earlier version of this  paper was presented at the
1st IEEE  Australian  Control Conference,  Melbourne, 2011 \cite{SNP}.}   }
\date{ \ }
\maketitle

\vspace{-1.5cm}
\begin{abstract}
 We consider the classic problem of pole placement by  state feedback. We adapt  the Moore  eigenstructure assignment algorithm  to obtain  a  novel parametric  form  for the  pole-placing gain matrix, and introduce an  unconstrained  nonlinear  optimization  algorithm to obtain a gain matrix   that will deliver robust  pole  placement.  Numerical experiments indicate the algorithm's performance compares favorably  against several  other notable robust pole placement methods from the literature.
\end{abstract}

\vspace{-.5cm}
\section{Introduction}

\label{secintro}

We consider the classic problem  of  pole placement for  LTI systems in state space form
  \bea
 \begin{array}{lcr}
 \dot x(t) = A\,x(t)+B\,u(t),\;\quad \label{syseq1}
  \end{array}
 \label{sys}
 \eea
where, for all $t \in \real $, $x(t) \in \real^n$ is the state, and  $u(t) \in \real^m$ is the control input.  $A$ and   $B$ are appropriate dimensional constant matrices. We assume that  $B$ has full column rank.  We let $\mathcal{L} = \{\lambda_1, \ldots, \lambda_\nu \}$ be a self-conjugate set of $n$ complex numbers,  with  associated  algebraic multiplicities $\mathcal{M} = \{m_1, \dots, m_\nu\}$  satisfying $m_1 + \dots + m_\nu= n$.
%  and $m_i \leq m$ for all $i \in \{1,\dots,\nu\}$.
  The  problem of  \textit{exact pole placement by state feedback} (EPP) is that of finding a real matrix $F$ such that  the closed-loop matrix $A + BF$ has non-defective eigenvalues in $\mathcal L$, i.e $F$ satisfies
\be \label{eppeq}
(A+BF)X =  X\Lambda
\ee
where $\Lambda$ is a $n \times  n$ diagonal matrix obtained from the  eigenvalues of  $\mathcal{L}$, including multiplicities,  and  $X$ is a non-singular  matrix of closed-loop eigenvectors of unit length. If $(A,B)$ has any  uncontrollable  modes, these are assumed to be  included within the set $\mathcal{L}$.
 The EPP problem has  been studied for several  decades,  and the existence of such a matrix yielding diagonal $\Lambda$  requires the  $m_i$ to  satisfy  certain inequalities in  terms of the controllability  indices of the  pair $(A,B)$ \cite{Rosenbrock-70}; in particular $m_i \leq m$ for all $m_i \in \mathcal M$ is required.  In  this paper we shall  assume  $(A,B,\mathcal L,\mathcal M)$ are such that at least one $F$ exists that  yields diagonal $\Lambda$.  Notable early papers  offering  algorithms to obtaining the required gain matrix $F$  include \cite{A},  which  gave a  method for  single-input single output (SISO) system, but this  was often  found to be numerically  inaccurate. Varga \cite{V1} gave a  numerically reliable method to obtain  $F$ for multiple-input multiple output (MIMO) systems.

 For SISO systems, $F$ is unique, while  for MIMO  systems it is not, and this naturally invites the selection  of $F$ that achieves the desired pole placement and also possesses other desirable  characteristics, such as minimizing the control input amplitude used, and improving numerical stability.  In order to  consider optimal  selections for  the gain matrix, it  is  important to have a  parametric formula  for the  set of gain matrices   that deliver the desired  pole  placement, and  numerous   such  parameterizations have appeared.  Bhattacharyya  and de Souza  \cite{BS} gave a  procedure for obtaining the gain matrix by solving a Sylvester equation in terms  of a $n \times  m $ parameter  matrix,  provided the  closed-loop  eigenvalues did not coincide with the open loop ones. Fahmy and O'Reilly   \cite{FO83}, gave a parametric form  in  terms  of the inverses of the  matrices $A -\lambda_i I$, which also required  the assumption  that the closed  loop  eigenvalues were all distinct from  the  open loop ones.  Kautsky \textit{et al} \cite{KNV} gave  a parametric form  involving  a QR-factorization  for  $B$  and a Sylvester equation for  $X$; this  formulation  did not require the  closed-loop poles to  be different from the  open-loop poles.

The classic eigenstructure assignment algorithm of  B.C. Moore \cite{M} quantified the freedom to simultaneously assign both the closed-loop eigenvalues, and also  select the associated eigenvectors.  As such it implicitly solved the EPP problem, but it did not   explicitly provide a parametric  formula   for the pole-placing  matrix,  nor did it  address any  optimal pole  placement problem. In this paper we adapt Moore's algorithm to obtain a simple parametric formula for the  pole-placing gain matrix, in  terms  of an $n \times  m $ parameter  matrix. The method obtains the  eigenvector matrix  $X$  by selecting eigenvectors from the nullspaces of the  system  matrices,  and thus avoids the  need for coordinate transformations.

The \textit{robust exact pole placement problem} (REPP) involves solving the EPP problem and also  obtaining $F$ that renders  the eigenvalues of $A + BF$  as insensitive to perturbations in $A$, $ B$ and $F $ as possible.
 Numerous results \cite{Wa} have appeared linking the sensitivity of the eigenvalues to various measures of the conditioning  of $X$, in terms of the Euclidean and Frobenius norms.  This classic optimal  control problem  also  has an  extensive  literature,  and typically  two approaches have  been used to  obtain good robust  conditioning.

Perhaps the best-known method for the  REPP is that  of  Kautsky  \textit{et al}  \cite{KNV}, which involved  selecting   an initial  candidate set of closed-loop  eigenvectors and then using a  variety of heuristic methods to make these vectors more orthonormal. This method   has been implemented as  MATLAB$^{\textrm{\tiny{\textregistered}}}$'s  \textit{place} command; this  implementation includes a heuristic extension to accommodate complex conjugate pairs in $\mathcal L$. This algorithm is also the basis of MATHEMATICA$^{\textrm{\tiny{\textregistered}}}$'s  \textit{KNVD} command. The  use of the  \textit{place}  algorithm has become  wide-spread in the control systems literature, and  introductory  texts advocating its use include  \cite{FPE} and  \cite{SSSH}, among many  others.

Since the  publication of  \cite{KNV}, many alternative  methods have  been proposed for the REPP.
 Tits and  Yang  \cite{TY} revisited the  heuristic   methods  of  \cite{KNV} and  offered  a  range  of improvements; the  algorithms were shown to be  globally  convergent.   Byers and Nash  \cite{BN},  Tam and Lam \cite{TL} and  Varga \cite{V2} cast the  problem  as an unconstrained nonlinear optimization problem, in  terms of the  Frobenius  conditioning,  to  be solved by  gradient iterative  search methods.  \cite{Ch} introduced a method for  minimizing the 'departure from normality' robustness measure, which considers the  size of the upper triangular part of the Schur form. Ait Rami \textit{et al}  \cite{RFBT} introduced a global constrained nonlinear optimal  problem in  terms of a Sylvester equation and showed that the solution  could be approximated by a convex linear problem  for which the authors gave an LMI-based algorithm.

%  It is well-known that  pole assignment problems are generally non-convex and  hence may possess many local minima. The heuristic approaches are known to computationally efficient, but cannot be expected to yield globally or even locally optimal  solutions. While the gradient iterative methods will obtain  a local minimum,  different  initial  conditions can yield  greatly  differing local  minima. The use of many  different initial conditions can lead to close estimates of the global minimum.

Various authors have  provided surveys comparing the performance of several of these algorithms. Sima \textit{et al} \cite{STY} conducted testing of  the algorithms from \cite{V1}, \cite{KNV} and \cite{TY} on  collections of systems  of varying   dimensions;  they concluded that the method of \cite{TY} generally gave superior Euclidean (2-norm) conditioning and also improved accuracy. \cite{Ch} considered the eleven benchmark  systems in the Byers-Nash collection (see Section \ref{seccomp} for a discussion of this collection), and compared the  author's  proposed methods, based on the Schur form of the open loop  systems, with those of  \cite{KNV} and \cite{TY} against a range of robustness measures. The methods of \cite{Ch} generally gave inferior results to those of
\cite{KNV} and \cite{TY}, with respect to the Frobenius conditioning. \cite{RFBT} tabulated figures results for the Frobenius conditioning performance of methods \cite{KNV}, \cite{TY}, \cite{BN} and  \cite{V2}. However, the  conditioning values were compiled directly from these  papers. Since some of these methods  were introduced into the  literature  more than  two  decades ago,  and noting that computational resources have improved dramatically over this time, using values from  original publications may unfairly disadvantage the earlier methods, in particular  \cite{BN}.

In this  paper we add to  this  extensive literature in  several  ways. In Section 2 we introduce our parametric form for the pole-placing gain matrix that solves the EPP. The formula is an adaptation of the pole placement method of Moore \cite{M}; the novelty here is  to use   Moore's  method  to obtain  a  parametric formula for  both   $X$,  the  matrix of   eigenvectors and   $F$,  the  pole-placing  gain matrix. We further show the  parametric form   is comprehensive, in  that it generates all possible $X$ and $F$ that solve (\ref{eppeq}), for the case where the eigenvalues have  multiplicity  of  at most $m$. In Section 3 we utilize this  parametric form to propose an  unconstrained optimization problem to  seek solutions to the REPP, to be solved by gradient search methods. Our approach most closely resembles that of  \cite{BN}, but with a different parametric formulation for the pole-placing gain matrix.

In Section 4 we select  five of  the most prominent methods for the REPP \cite{KNV}, \cite{TY},  \cite{V2}, \cite{BN}  and \cite{RFBT},  and conduct extensive numerical testing  to compare their  performance  against  our method. The first three of these were chosen as they are widely  used in  the forms of the  \matlab  toolboxes \textit{place}, \textit{robpole} and \textit{sylvplace} respectively.  \cite{BN} has attracted  a large  number of citations over more than two  decades, and \cite{RFBT} is the  most recent  publication to offer a novel approach for the  REPP.  All methods were  implemented in \matlab 2012a, running  on the same computing  platform. In addition to  conditioning, we also compare their accuracy, matrix gain and runtime.  Finally, Section 5 offers some conclusions  as to  the  relative performance of these six  methods; our   method  will be  shown  to offer some performance advantages over all the other methods surveyed.

\section{Pole placement via Moore's algorithm}
\label{Moore}

We now revisit Moore's method \cite{M} and adapt it to  give  a  simple parametric  formula  for a gain matrix $F$ that solves the  pole placement problem, in terms  of an  arbitrary   real parameter matrix.  We begin with some  definitions and notation.  For each  $i \in \{1,\dots,\nu\}$,  we define the $n \times (n +m)$ system matrix
 \be \label{Slambda}
 S(\lambda_i)  =  [ A-\lambda_i I_n \ B]
 \ee
 where $I_n$ is  the identity matrix of size $n$.  We let $T_i$ be a basis matrix for the  nullspace of  $S(\lambda_i)$,  we  use $s_i$ to denote the dimension of this  nullspace,  and we denote $T= : [T_1 \dots T_\nu]$. It follows that $s_i = m$, unless  $\lambda_i$ is an uncontrollable mode of the  pair $(A,B)$, in which case we will have  $s_i > m$.  Let $M$ denote any complex  matrix partitioned into submatrices $M= [ M_1 | \dots | M_\nu]$ such that any complex submatrices occur consecutively in complex conjugate pairs. We define a  real matrix $Re(M)$ of the same dimension  as  $M$ thus: if $M_i$ and $M_{i+1}$ are consecutive complex conjugate submatrices  of $M$, then the corresponding submatrices  of  $Re(M)$ are $\frac{1}{2}(M_i+M_{i+1})$ and $\frac{1}{2j}(M_i-M_{i+1})$.
    Finally, for any real or complex  matrix $X$ of with at  least $n+m$ rows, we define matrices $\overp(X)$ and  $\undp(X)$  by taking the first $n$ and last $m$ rows of $X$, respectively.

 \begin{prop} \label{Mooreprop}
 Let the eigenvalues  $\{\lambda_1,\ldots,\lambda_\nu\}$ be ordered so that,  for some  integer $s$, the first $2s$ values are complex while the remaining are real,  and for all odd $i \le 2\,s$ we have $\lambda_{i+1}=\bar \lambda_i$.   Let $ K:=\diag(K_1, \dots, K_\nu)$,
 %be  a real  matrix of  dimensions $m \times n$  partitioned so that
 %\be \label{Keq}
% K = [k_1 | \dots | k_\nu],
% \ee
  where  each $K_i$ is of dimension $s_i \times m_i$, and for all odd $i \leq 2s$, we have $K_i=\bar K_{i+1}$.
 Let $M(K)$ be  an $(n+m) \times n$  complex matrix  given  by
 \be \label{Meq}
 M(K) =    TK
 \ee
 and let
 \bq \label{XVWeq}
  X(K) &= & \overp(M(K)),  \label{Xeq}    \\
  V(K) &= &  \overp(Re(M(K))) \label{Veq} \\
  W(K) &= &  \undp(Re(M(K))) \label{Weq}
 \eq For almost every choice of the parameter matrix $K$, the rank of  $X$ is equal to $n$. The  set of all $m \times n$  gain matrices $F$ satisfying (\ref{eppeq}) is parameterised in $K$ as
 \be \label{Feq}
 F(K) = W(K)V(K)^{-1}
 \ee
where $K$ is such that  $\rank(X(K))=n$.
\end{prop}

 \IEEEproof  For any  given $K$, let $M(K)$ be  partitioned according to
 \be
  M(K)  = \left[ \begin{array}{ccc} V_1^\prime & \dots & V_\nu^\prime \\  W_1^\prime & \dots & W_\nu^\prime \end{array} \right]
  \ee
  where  each $V_i^\prime $ and  $ W_i^\prime $ are  matrices of  dimensions $n \times m_i$ and  $m \times m_i$
  respectively,   such that
\be
\label{ML}
(A -\lambda_i\,I_n) V_i^\prime +  B W_i^\prime =0
\ee
Note that, for odd $i \leq  2s$, we have that $V_i^\prime  = \bar V_{i+1}^\prime$ are conjugate matrices, as  $K_i = \bar K_{i+1}$. Moreover,  since  $\mathcal L$ is symmetric, we also  have $m_i = m_{i+1}$. Define real matrices
\bea
\label{defV}
V_i=\left\{ \begin{array}{lll}
\frac{1}{2} (V_i^\prime+V_{i+1}^\prime) & \textrm{if $i \le 2\,s$ is odd}, \\
\frac{1}{2j} \,(V_{i-1}^\prime-V_{i}^\prime) & \textrm{if $i \le 2\,s$ is even}, \\
V_i^\prime & i >2\,s
\end{array} \right.
\eea
and define $W_i$ similarly. Then matrices $X$, $V$ and $W$ in (\ref{Xeq})-(\ref{Weq}) may  be written as

\noindent $X = [ \,V_1^\prime \;\; V_2^\prime\;\; \ldots\;\; V_{2\,s}^\prime\,| \,V_{2\,s+1}^\prime \;\; V_{2\,s+2}^\prime\;\; \ldots\;\; V_{\nu}^\prime\,]$,
 $V=[\,V_1 \;\; V_2\;\; \ldots\;\; V_{2\,s}\,| \,V_{2\,s+1} \;\; V_{2\,s+2}\;\; \ldots\;\; V_{\nu}\,]$.

 \noindent  and   $W=[\,W_1 \;\; W_2\;\; \ldots\;\; W_{2\,s}\,| \,W_{2\,s+1} \;\; W_{2\,s+2}\;\; \ldots\;\; W_{\nu}\,]$.  Let
 \be \label{Peq}
 R_i = \frac{1}{2}  \left[ \begin{array}{cc} I_{m_i} & -jI_{m_i} \\  I_{m_i}  & jI_{m_i} \end{array} \right]
 \ee
 Then for each odd $i \leq  2s$, we have $[ V_i^\prime \;\; V_{i+1}^\prime ] R_i = [ V_i \;\; V_{i+1} ]$ and  $[ W_i^\prime \;\; W_{i+1}^\prime ]R_i = [ W_i \;\; W_{i+1} ]$. Now assume $K$ is such  that $\rank(X(K))=n$; then $V(K)$ is  non-singular, and we can obtain  $F$ in  (\ref{Feq}).
 %To  see that $F$   satisfies  (\ref{eppeq}),   we
   We obtain  $F[V_i^\prime  \ V_{i+1}^\prime ] = [W_i^\prime  \ W_{i+1}^\prime ]$ for  odd $i\in\{1,\ldots,2\,s\}$  and    $FV_i^\prime  = W_i^\prime$ for all $i \in  \{2s+1, \dots, \nu\}$. Hence  (\ref{ML}) can be written as
\bq
 (A+BF)\left[\begin{array}{cc} V_i^\prime & V_{i+1}^\prime \end{array} \right]
  & =  & \left[ \begin{array}{cc}   V_i^\prime & V_{i+1}^\prime \end{array} \right]\diag(\lambda_i I_{m_i},  \lambda_{i+1} I_{m_{i}}),
   \mbox{for odd $i\in\{1,\ldots,2\,s\}$}   \label{clcomplexi}   \\
(A+B\,F)V_i^\prime & = &  V_i^\prime(\lambda_i I_{m_i}),  \mbox{for  $i\in\{2\,s+1,\ldots,\nu\}$}, \label{clreali}
    \eq
Thus we obtain (\ref{eppeq}). To  see that  this  formula is comprehensive, we let $F$  be  any  real  gain matrix  satisfying  (\ref{eppeq}). The nonsingular  eigenvector  matrix $X$ is  comprised of column  vectors  $V_i^\prime$ of dimension $n \times m_i$ corresponding to each  eigenvalue,  such  that (\ref{clcomplexi}) and (\ref{clreali})  hold. Applying  $F[V_i^\prime  \ V_{i+1}^\prime ] = [W_i^\prime  \ W_{i+1}^\prime ]$ for  odd $i\in\{1,\ldots,2\,s\}$  and $FV_i = W_i$ for all $i \in  \{2s+1, \dots, \nu\}$,   we  obtain $V_i^\prime$ and $W_i^\prime$ such  that (\ref{ML})  holds. Thus each  column  vector of the   matrix $[V_i^\prime \ W_i^\prime]^T$ lies in the  kernel of  $S(\lambda_i)$,  and we have a  coefficient  vector $K_i$ such  that  $[V_i^\prime \ W_i^\prime]^T = T_iK_i$. The  complex conjugacy of   $V_i^\prime$ and $ V_{i+1}^\prime$,  for each  odd $i\in\{1,\ldots,2\,s\}$, implies the  conjugacy of  $K_i$ and $ K_{i+1}$.  Thus
   we obtain $M(K)$ in (\ref{Meq}) yielding  $F$ in  (\ref{Feq}).

Finally we let  $K$ be arbitrary parameter matrix and consider the rank of  $X(K)$.  We introduce  $\Phi=\overp(T)$  and denote  $\Phi_{1},  \ldots, \Phi_{\nu}$ as a basis for $\im \Phi$.
  If $\rank (X(K))$ is smaller than $n$, then one column of the matrix  $[\Phi_{1} K_{1,1}  \dots \Phi_\nu K_{\nu,m_\nu}]$  is linearly dependent of all the remaining ones. (Here we have  used  $K_{i,j}$ to  denote the  $j$-th column of $K_i$).   For brevity, let us assume this is the last column. Then there exist $n-1$ coefficients $\alpha_{1,1},\ldots,\alpha_{\nu,m_\nu-1}$ (not all equal to zero) for which
  \be  \label{Knu}
  \Phi_{\nu}  K_{\nu,m_\nu}
    = \sum_{i=1}^{\nu -1 } \sum_{j=1}^{m_i} \alpha_{i,j}\,\Phi_{i}K_{i,j} + \sum_{j=1}^{m_\nu -1} \alpha_{\nu,j}\,\Phi_{\nu}K_{\nu,j}
  \ee
has a unique solution in $K_{\nu,m_\nu}$. As $K_{\nu,m_\nu}$ is an $s_\nu$-dimensional parameter vector, (\ref{Knu}) constrains $K_{\nu,m_\nu}$ to lie upon  an  $(s_\nu -1)$-dimensional hyperplane, which has empty interior.  Thus the  set of parameters $K$  that lead to a loss of rank in $X(K)$ is given by the union of at  most $n$ hyperplanes of empty interior. This set therefore has empty interior, and thus also zero Lebesgue measure. Thus we see that $X(K)$ and   hence  $V(K)$ are non-singular for almost all choices of the parameter matrix $K$.
 \endproof

 The above formulation takes its inspiration from  the  proof of Proposition 1  in \cite{M}, and hence we  shall  refer to  (\ref{Xeq})-(\ref{Feq}) as the  {\it Moore parametric form}   for  $X$  and  $F$. We note however  that \cite{M} only considered the case of distinct  eigenvalues, and did not offer any explicit parametric formula for  the  pole-placing gain matrix. Moreover, it did not  show that  all  matrices $X$ and $F$  solving   (\ref{eppeq})  could  be parameterized in the above manner.

It is  interesting to compare this parametric form with that of  \cite{KNV}, in  which the  eigenvectors  comprising   $X$ were  obtained from the nullspaces of the matrices  $U_1(A-\lambda_i I)$, where  the  parameter $U_1$ was obtained  from  the QR-factorization for $B = [U_0 \ U_1][ Z \  0]^T$, and was also  required to  satisfy $U_1(AX - X\Lambda) =0 $. By  contrast,  the Moore parametric  form  obtains the eigenvectors  directly from the nullspaces  of the system matrices  $[A- \lambda_i I_n  \ B]$.

\section{Robust and  minimum gain pole placement}

 When  $A + BF$ has $n$ distinct eigenvalues, the sensitivity of an eigenvalue $\lambda_i$ of $A+BF$ to perturbations in $A$, $B$, and $F$ can be represented by the condition number \cite{Wa}
\be\label{ci}
    c_i=\frac{\|y_i\|_2\|x_i\|_2}{|y_i^Tx_i|}
\ee
where $y_i$ and $x_i$ are the left and right  eigenvectors associated with $\lambda_i$.  For the case where  $A + BF$ is  non-defective but  has  repeated eigenvalues, see \cite{Sun} for a definition of the corresponding condition numbers. Furthermore, we have \cite{KNV}
\be
    c_\infty:=\max_i{c_i}\leq \kappa_2 (X)\leq\kappa_{fro}(X)
\ee
where $\kappa_2 (X)= \|X\|_2\|X^{-1}\|_2$ and $\kappa_{fro} (X)= \|X\|_{fro}\|X^{-1}\|_{fro}$ are the condition numbers of the  matrix of eigenvectors $X$ with respect to the Euclidean and Frobenius norms. Following  \cite{RFBT,BN,TL}, we propose to  address  the  REPP  problem by  minimizing the condition number of $X$ with respect to Frobenius norm.  The  objective function to  be  minimized is
\be
    f_1(K)=\kappa_{fro} (X(K))=\|X(K)\|_{fro}\|X^{-1}(K)\|_{fro}
\ee
where the input parameter matrix  $K$ is defined as in Proposition \ref{Mooreprop}. Note it is possible to reduce the Frobenius norm  of a matrix $X$ by suitably scaling the  lengths of its  column vectors. When $X$ is the solution to (\ref{eppeq}), such  scaling does not improve the eigenvalue conditioning  in (\ref{ci}). Hence we assume that the column vectors of $X$  have been normalised.

As pointed out in \cite{BN}, for efficient computation we can study an alternative objective function
\be
    f_2(K)=\|X(K)\|_{fro}^2+\|X^{-1}(K)\|_{fro}^2
\ee
because the two objective functions are equivalent. An imported related problem is that of minimizing the norm of the gain matrix $F$. The \textit{minimum gain robust exact pole placement problem} (MGREPP) involves simultaneously minimizing both the conditioning and the matrix gain  via the weighted  objective function
 \be
    f_3(K)=\alpha \kappa_{fro}(X(K)) +(1 - \alpha)\|F(K)\|_{fro} \label{mgrepp}
\ee
where $\alpha$ is a weighting factor, with $0 \leq \alpha \leq 1$.
 %For $\alpha  = 0$ this defines a pure matrix gain minimization problem, while for $\alpha =1$ we have the REPP problem.
   Minimizing $f_3$ involves a gradient search employing the first and second order derivatives of $ \kappa_{fro}(X(K)) $ and $\|F(K)\|_{fro}$; expressions for these were given in \cite{SNP}.

\section{Performance comparison of robust pole placement methods} \label{seccomp}

In this section we conduct  extensive numerical experiments to  compare the  performance  of our method against those of \cite{KNV}, \cite{RFBT}, \cite{TY}, \cite{BN} and \cite{V2}. To provide a comprehensive contemporary  survey, we implemented these algorithms on the same modern computer, an Intel$^\registered$ Core$^\trademark$  Quad CPU, Model Q9400 at  2.66 GHz with   3326 MB of RAM running Windows$^\trademark$ XP  and MATLAB$^{\textrm{\tiny{\textregistered}}}$ 2012a. Implementation of  \cite{KNV}   was done with MATLAB$^{\textrm{\tiny{\textregistered}}}$'s  \textit{place} command. For \cite{TY} and  \cite{V2}, we used the \textit{robpole} and \textit{sylvplace} MATLAB$^{\textrm{\tiny{\textregistered}}}$ toolboxes, kindly  provided to us by the authors.  For \cite{BN},   \cite{RFBT} and  our  own method, we wrote MATLAB$^{\textrm{\tiny{\textregistered}}}$ toolbox implementations for each. The \cite{RFBT} algorithm  requires an  LMI solver; we chose the public-domain \textit{cvx} toolbox \cite{GB}. We shall  refer to these as \textit{byersnash}, \textit{rfbt} and $\Span$ (our  own method). The names are derived from the  names of the respective authors.

To obtain a fair comparison between these methods, we need to consider the runtime allocated to them. The methods of \cite{BN},  \cite{V2} and  our proposed  method all employ gradient iterative searches, so the values they deliver are contingent upon the initial condition (input parameter matrix $K$) used. The $\sylv$  toolbox randomly generates an  initial condition, and thus offers different outputs (different $F$) each time it is run. To obtain repeatable results, we provided the $\byers$ and $\Span$ toolboxes with a pre-specified collection of  input parameter matrices $K$ composed of canonical vectors. The output shown from  each  of \textit{byersnash}, \textit{sylvplace} and $\Span$  is the best  result  from  all the  initial conditions searched within the allocated  runtime. By contrast   $\place$, $\rob$ and  $\rfbt$ all employ a designated starting point, and hence their runtime is simply the time taken to execute their method.

 \subsection{Robust conditioning comparison  using the Byers and Nash benchmark examples}

 Byers and  Nash  \cite{BN} gave  a collection of eleven  benchmark  example systems, and many authors, including  \cite{TY}, \cite{V2} and \cite{RFBT} used these examples to compare the performance of their pole placement  methods. Following this  well-established tradition, our first set  of  comparisons employs  these well-known  examples.   The results   are given in  Table \ref{Surv1}.
  We have  used $\kappa_{fro}(X)$ as  the  performance measure, and we also  show the matrix  gain  used.

  The average runtimes for $\place$, $\rob$  and  $\rfbt$  for the 11 sample systems were $ 0.05$,  $0.095 $ and $ 14.1$ seconds, respectively. For  $\byers$, $\sylv$  and  $\Span$ we arbitrarily set the runtime to be $n$ seconds, where $n$ is the system dimension, leading to  average runtimes of 4.5 seconds, this being  the average of the system dimensions in  the collection.

 Ignoring differences  in the conditioning of smaller than 1\%, we  conclude that \textit{byersnash} and \Span\ had the  best  or equal  best conditioning in  all 11 examples.  \textit{sylvplace} and  \textit{rfbt} had the best or equal  best in 7 cases, while \textit{robpole} had  best  or equal best in  5 cases. Finally  \textit{place} gave the best or equal best in  4 cases.  $\place$ and $\rob$  had the shortest runtimes, while  \textit{rfbt} had noticeably the longest.   We note that the conditioning numbers given here differ significantly from those that were  published in \cite{BN} and \cite{RFBT}. This  may be explained by the fact that these authors did not  require  the columns of  $X$ to be of  unit length. Since methods \cite{KNV} and \cite{TY} normalise the columns of $X$, this is essential for a fair comparison of all six methods.

 \subsection{Robust conditioning comparison with sets of higher-dimensional systems}

 To probe more deeply into the performance delivered by these six methods, we need to move beyond the low-dimensional  examples in  the  Byers and Nash collection. In Survey 2 we generated three  sets of 500 sample systems with $(A,B)$,  all of  state dimension $n =20$, and with  control input dimensions of  $m=2$, $m=4$ and $m=8$. The  pole positions $\mathcal L$ were chosen to be all distinct, with  a mixture of  real and complex values. The entries of $A$,  $B$ and $\mathcal L$ took uniformly distributed   values  within the  interval $[-2,\ 2]$.
To  compare the conditioning, accuracy, and matrix gain of each method, we computed,  for each system $j \in \{1, \dots,500\}$ and each  method  $\star \in \{\textit{place},\textit{robpole}, \textit{byersnash}, \textit{sylvplace}, \textit{rfbt}, \textit{span}\}$,
\bi
\item $\kappa_{fro}(\star,j)$: the Frobenius conditioning of method $\star$  for the  $j$-th system;
\item $c_\infty(\star,j)$: the $c_\infty$  conditioning of method $\star$  for the  $j$-th system;
     \item $\Delta(\star,j)$: the accuracy of method $\star$ on the $j$-th  system, equal to
      %$\max\{|eig_i(A+BF)-\lambda_i| :  i \in \{1,\dots,n\}\} $,
        the  largest absolute  value difference  between  each eigenvalue of  $A + BF$ and the corresponding $\lambda_i$  in  $\mathcal{L}$.
     \item $\|F\|_{fro}(\star,j)$: the Frobenius norm of $F$ from Method $\star$ on system $j$.
     \ei
  Noting that  \textit{place} is the industry standard for the  REPP, we chose to compare all the other methods according to their ability to improve upon \textit{place},  and computed  comparative performance indices relative to  \textit{place} for each  method,  and for each performance criterion, as follows:
\bq
   (1 - \ind(\star,\kappa_{fro}))^{500}  &  =  &  \prod_{j=1}^{500} \ \frac{\kappa_{fro}(\star,j)}{\kappa_{fro}(place,j)} \label{indkappafro} \\
   (1 - \ind(\star,c_\infty))^{500} & = & \prod_{j=1}^{500} \ \frac{c_\infty(\star,j)}{c_\infty(place,j)}  \\
   (1 - \ind(\star,\Delta))^{500}    &   = &  \prod_{j=1}^{500} \ \frac{\Delta(\star,j)}{\Delta(place,j)} \\
   (1 - \ind(\star,\|F\|_{fro}))^{500} & =  & \prod_{j=1}^{500} \ \frac{\|F\|_{fro}(\star,j)}{\|F\|_{fro}(place,j)} \label{indFfro}
\eq

For example, in (\ref{indFfro}),   if $\ind(\rob,\|F\|_{fro}) = 0.1$,  then Method $\rob$ gives values of $ \|F\|_{fro}$ that are on average 10\% smaller  than  \textit{place}.  Larger indices imply greater improvement on \textit{place},  and negative indices indicate performance inferior to \textit{place}. The local gradient search  methods $\Span$, \textit{byersnash} and  \textit{sylvplace} were each given $20$  seconds of runtime per sample system;  the results shown in Table \ref{Surv2} represent the best conditioning  performance achieved from all  the initial  conditions searched within that time  period. For $\rob$  and  $\rfbt$, the average runtime  per sample system were $ 0.552$ and $ 125$ seconds ($m=2$), $0.552 $  and $ 82.9$  seconds ($m=4$),
 and $0.552$  and  $ 55.2$ seconds ($m=8$).

 The results  show that the best performance for robustness and   gain minimisation  were given by  $\Span$,  $\byers$ and $\sylv$. Both $\sylv$ and  $\rfbt$ were less accurate than $\place$, by  several orders  of  magnitude in  the case  of   $\rfbt$, which also required  substantially  longer runtime. While  all methods offered improved conditioning  with reduced gain over $\place$,  this was reduced for the larger values  of $m$,  which may  be attributed to the improved performance of  $\place$ when it has more  control inputs to work with.

\subsection{Weighted gain minimisation and conditioning problem}

  Among the methods in our survey, only \cite{V2} ($\sylv$) considered  the MGREPP problem (\ref{mgrepp}). Our   Survey 3 compares the performance of  $\sylv$ and $\Span$  for the same 500 sample systems used in Survey 2,  with $m=2$,  for several  different values of the weighting factor  $\alpha$.  We  again gave  $\Span$ and $\sylv$ 20 seconds of  runtime per sample system,  and  computed the performance improvement indices (\ref{indkappafro})-(\ref{indFfro}) relative to the gain matrix  delivered by $\place$;  again larger figures indicate greater improvement. The results are shown in Table \ref{Surv3}. Both methods were able to offer significant reductions in gain, at the  price of some reduction in the robustness measures, relative  to  the pure robustness  problem ($\alpha =1$). However  $\Span$ did  so  with  far superior accuracy.
Considering  the impact of different values of the weighting factor,  we see that for $\alpha = 0.1$,  there was little difference in the conditioning, and only slight improvement in the  matrix  gain.  For $\alpha \ra 0$ we observed up  considerable reduction in the matrix gain, but this eventually comes at the cost of significantly inferior conditioning. These results suggest values around $\alpha = 0.001 $ can  give a good balance between these two criteria.

   \subsection{Systems with uncontrollable modes}

 The  EPP  problem remains well-posed   for systems with  uncontrollable  modes, provided these are included within the  set $\mathcal L$. The  methods  $\place$,  $\sylv$,  $\rob$, $\rfbt$ all  assumed controllability of the system, as  part of their  problem formulation. In principle this involves no loss of generality, since  the application of a Householder staircase transformation  can decompose any system into its  controllable and  uncontrollable  parts.  Nonetheless is it is interesting to  consider the  ability of these toolboxes to  accommodate  uncontrollable  modes. In our  final survey, we obtained  100 systems  $(A,B)$,  with $n=3$ and  $m =2$, that  contained  one  uncontrollable mode. We then  chose $\mathcal L$ to include this  mode, plus one  pair of complex conjugate  modes. We defined failure to solve the EPP as being any one  of (i) an error was returned upon execution of the algorithm, (ii) any of the closed-loop poles differed by more than $5\%$ from their  desired location,  and (iii) the gain of $F$ was undefined or greater than $10^{10}$.  We observed failures as  follows: $\place$,  $\sylv$,  $\rob$  and  $\rfbt$ had 100, 98,  30 and  12 failures,  respectively;  we concluded these toolboxes in their present form cannot reliably solve the EPP in these conditions.   $\byers$ and our method $\Span$ had  no   failures; we attribute their  superior reliability  to their usage of  nullspace methods.  Uncontrollable modes increase the  column dimension of the corresponding nullspace basis matrix; for  $\byers$ and  $\Span$ this is readily accommodated by adjusting the row dimension of the parameter matrix.

\section{Conclusion}
 We have introduced a parametric formula  for the exact pole placement of linear systems via  state feedback, derived from Moore's classic eigenstructure method. This parametric form was used to formulate the robust and  minimum  gain exact pole placement problem  as an unconstrained  optimization problem, to  be solved by gradient iterative methods.

 The  method was implemented as  a MATLAB$^{\textrm{\tiny{\textregistered}}}$ toolbox called $\Span$, and  its  performance  was compared against several other methods from the classic and  recent literature. All methods considered gave superior performance to  the widely  used MATLAB$^{\textrm{\tiny{\textregistered}}}$ \texttt{place} command, albeit with somewhat longer runtime.  When  the  Frobenius conditioning of the eigenvector  matrix is used as the  robustness measure, the  best performance was provided by the our proposed method,  and also the Byers-Nash method. The results suggest that, in  comparison  with heuristic methods,  gradient iterative methods are best able to take advantage of the high levels of computational power that are now widely available. They also suggest that methods based on nullspaces of appropriate system matrices may offer superior accuracy of pole  placement to  those adopting Sylvester matrix  transformations.

 For  a given system $(A,B,\mathcal L,\mathcal M)$, $\byers$ and $\Span$  will in general yield quite different gain matrices, offering different performance values, so both methods should be considered for optimal performance. While Byers  and Nash considered only  the  robustness,   our method is able to accommodate  a  combined robustness and gain minimization  approach, enabling the designer to obtain  significantly reduced  gain in  exchange for somewhat inferior conditioning.

 The authors would like to thank Andre Tits and Andreas  Varga for  providing us with  copies of their \textit{robpole} and  \textit{sylvplace}  toolboxes,  and Ben Chen  for bringing the  classic eigenstructure assignment paper  by B.C. Moore \cite{M} to our attention. We also  thank the anonymous reviewers for some constructive suggestions.

\bibliographystyle{plainnat}

\begin{table}[tbh!]
\begin{center}
\caption{Survey 1: REPP with the Byers  Nash Examples} \label{Surv1}

\begin{tabular}{||c||c|c||c|c||c|c||}
  \hline
  Example &  \multicolumn{2}{c||}{$\place$\cite{KNV}}  &  \multicolumn{2}{c||}{$\byers$\cite{BN}} &  \multicolumn{2}{c||}{$\rob$ \cite{TY}}
   \\ \cline{2-7}
    & $\kappa_{fro}(X)$ & $\|F\|_{fro}$ & $\kappa_{fro}(X)$ & $\|F\|_{fro}$ & $\kappa_{fro}(X)$ & $\|F\|_{fro}$ \\
     \hline
  1 & 6.5641 & 1.364 & 6.4451 & 1.4582 & 7.3214 &  1.3338 \\

  2 & 57.491 &  301.37 & 50.224 & 355.19 & 52.972 & 224.95 \\

  3 & 103.18 & 105.06 & 46.238 & 77.215 & 55.987 & 49.104 \\

  4 & 13.431 &  9.899 & 13.421 & 9.4485 & 13.421 & 9.4462 \\

  5 & 146.18 & 4.8496 & 142.39 & 4.5561 & 144.78 & 5.4168   \\

  6 & 6.0018 & 21.5 & 5.9633 & 23.25 & 6.0262 & 20.197  \\

  7 & 12.375 & 233.64 & 11.302 & 326.35 & 12.017 & 235.08   \\

  8 & 36.986 & 15.7600 & 6.1824 & 28.033 & 6.1824 & 28.599 \\

  9 & 28.682 & 2356.5 & 23.915 & 832.22 & 23.937 &  823.70 \\

  10 & 4.0029 & 1.4897 & 4.113 & 5.2687 & 4 & 1.5174 \\

  11 & 14618 & 6692.1 & 14510 & 6580.8 & 14510 & 6580.7 \\
\hline
  \hline
  Example &  \multicolumn{2}{c||}{$\sylv$\cite{V1}} &  \multicolumn{2}{c||}{$\rfbt$\cite{RFBT}} &  \multicolumn{2}{c||}{$\Span$}
   \\ \cline{2-7}
   & $\kappa_{fro}(X)$ & $\|F\|_{fro}$ & $\kappa_{fro}(X)$ & $\|F\|_{fro}$ & $\kappa_{fro}(X)$ & $\|F\|_{fro}$ \\
    \hline
  1 & 6.5997 & 1.4662 & 6.5595 & 1.5253 & 6.4451 & 1.4582 \\

  2 & 50.042 & 327.75 & 50.185 & 361.01 & 50.224 & 355.17    \\

  3 & 45.741 &  72.285 & 45.772 &  73.582 & 46.223 & 77.146 \\

  4 & 13.421 & 9.4465 & 13.421 & 9.366 & 13.421 & 9.4432   \\

  5 & 141.99 & 4.8472 & 142.82 & 4.3963 & 142.39 & 4.556    \\

  6 & 5.9361 & 22.474 & 6.4086 & 14.771 & 5.9622 & 23.318  \\

  7 & 11.353 & 271.17 & 12.280 & 297.85 & 11.301 & 271.06  \\

  8 & 6.1824 & 21.827 &   9.381 &  39.300 & 6.1824 & 21.102  \\

  9 & 24.23 & 903.11 & 23.925 & 884.84 & 23.916 & 831.23    \\

  10 & 4.113 & 1.513 & 4 & 1.5185 & 4 & 1.517  \\

  11 & 16571 & 10716 & 14475 & 6642 & 14510 & 6581.3 \\
\hline
\end{tabular}
\end{center}
\end{table}

\np

\begin{table}[tbh!]
\begin{center}
\caption{Survey 2: REPP with  higher-dimensional systems} \label{Surv2}
\begin{tabular}{||c|c||c|c|c|c|c||}
  \hline
 System Dimension & Metric & \multicolumn{1}{c|}{$\byers$\cite{BN}} & \multicolumn{1}{c|}{$\rob$ \cite{TY}} & \multicolumn{1}{c|}{$\sylv$\cite{V1}} & \multicolumn{1}{c|}{$\rfbt$\cite{RFBT}} &   \multicolumn{1}{c||}{$\Span$}
   \\ \cline{2-5}
     \hline
$n=20$,  & $\kappa_{fro}(X) (\%)$  & 54.670  & 9.8815 & 51.938 & 41.332 & 54.603 \\

 $m = 2$, & $c_{\infty} (\%)$  & 62.047 & 10.620 & 59.759 &  49.447 & 61.983  \\

 sys = 500  & $\|F\|_{fro} (\%)$ & 23.555 & 1.9292 & 22.310 & 14.337 & 23.276 \\

       &  Accuracy (\%) & 67.356 & 26.998 & -1.0082 & -46237 & 64.344  \\
\hline
  \hline
$n=20$,  &   $\kappa_{fro}(X) (\%)$ & 37.268 & 9.150 & 36.725 & 31.048 & 37.264 \\

 $m = 4$, & $c_{\infty} (\%)$  & 49.418 & 9.8601 & 50.226 & 43.374 & 49.400  \\

  sys =500  & $\|F\|_{fro} (\%)$ & 15.677 & 4.3745 & 15.524 & 11.163 & 15.698 \\

         & Accuracy (\%) & 45.057 & 23.760 & -65.586 & -169100 & 43.034 \\
\hline
  \hline
 $n=20$,  &   $\kappa_{fro}(X) (\%)$ & 15.198 & 7.7702 & 11.745 & 12.849 & 15.197 \\

  $ m = 8$, &  $c_{\infty} (\%)$  & 23.271 & 10.067 & 20.848 & 20.840 & 23.236  \\

  sys =500       & $\|F\|_{fro} (\%)$ & 3.7940 & 4.7471 & 3.3979 & 1.7034 & 3.7860 \\

         & Accuracy (\%) & 18.525 & 17.8859 & -44.635 & -338240 & 16.225 \\

\hline
\end{tabular}
\end{center}
\end{table}

\begin{table}[tbh!]
\begin{center}
\caption{Survey 3: MGREPP with higher dimensional  systems ($n = 20, m = 2$, sys =500)} \label{Surv3}
\begin{tabular}{||c||c|c|c|c|c|c|c||}
\cline{1-7}
\multirow{2}{*}{Metric}
      & \multicolumn{2}{c|}{$\alpha = 0.0001$}  &  \multicolumn{2}{c|}{$\alpha = 0.001$} &  \multicolumn{2}{c|}{$\alpha = 0.1$}
    \\ \cline{2-7}
    \cline{2-7}

  &  \multicolumn{1}{c|}{$\Span$}  &  \multicolumn{1}{c|}{$\sylv$\cite{V1}}& \multicolumn{1}{c|}{$\Span$} & \multicolumn{1}{c|}{$\sylv$\cite{V1}} &   \multicolumn{1}{c|}{$\Span$}  &  \multicolumn{1}{c|}{$\sylv$\cite{V1}}
   \\ \cline{2-7}
     \hline

  $\kappa_{fro}(X) (\%)$ & -25.578  & 23.980 & 37.641 & 41.906 & 53.936 & 51.699  \\

  $c_{\infty} (\%)$  & -13.540 & 33.929 & 45.966 &  51.379 & 61.213 & 59.465  \\

  $\|F\|_{fro} (\%)$ & 50.319 & 38.046 & 43.577 & 37.740 & 27.509 & 26.404   \\

  Accuracy (\%) & 16.992 & -46.326 & 57.833 & -16.025 & 65.643 & -1.0463 \\

\hline
\end{tabular}
\end{center}
\end{table}

\end{document}